\documentclass{svmult}
\usepackage{amssymb}
\usepackage{amsmath}
\usepackage{graphicx}

\title*{An introduction to arithmetic groups}
\author{Christophe Soul\'e}
\institute{CNRS and IHES, 35 Route de Chartres, 91440 Bures sur Yvette, France.  \texttt{soule@ihes.fr}}

\begin{document}

\maketitle

{\it Arithmetic groups} are groups of matrices with integral coefficients. 
They first appeared in the work of Gauss, Minkowski and others on the 
arithmetic theory of quadratic forms. Their {\it reduction theory} consists 
in showing that, after a linear change of variables with integral 
coefficients, any quadratic form can be forced
 to satisfy an appropriate set
of inequalities.

Around 1940, Siegel developed a general theory of arithmetic subgroups of 
classical groups, and the corresponding reduction theory. Later on, once 
Chevalley, Borel, Tits and others had developed the general theory of 
algebraic groups, one could speak of the arithmetic subgroups of any linear 
algebraic group  over ${\mathbb Q}$. Borel et al. extended the work of
Siegel to arbitrary arithmetic groups.

These groups play a fundamental role in number theory, and especially in the 
study of automorphic forms, which can be viewed as complex valued functions 
on a symmetric domain which are invariant under the action of an arithmetic 
group. In the last ten years, it appeared that some arithmetic groups are 
the symmetry groups of several string theories. This is probably why this 
survey fits into these proceedings.

In a first chapter we shall describe the classical reduction theory of 
quadratic forms. After describing the action of ${\rm SL}_2 ({\mathbb Z})$ 
on the Poincar\'e upper half-plane (Theorem 1) we explain how Siegel defined 
a fundamental domain for the action of ${\rm GL}_N ({\mathbb Z})$ on real 
quadratic forms in $N$ variables (Theorem 2). We then proceed with the 
general definition of linear algebraic groups over ${\mathbb Q}$ and their 
arithmetic subgroups (\S~3). An important example is a construction of 
Chevalley which defines an arithmetic group $G({\mathbb Z})$ when given any 
root system $\Phi$ together with
 a lattice between the root lattice and the weight
lattice of $\Phi$ (3.3). In 3.4 (and in the Appendix)
we explain how the group $E_7 ({\mathbb Z})$
of \cite{soule:HT} is an example of this construction. We then describe
the general construction of
Siegel sets and the reduction theory of arithmetic groups
(Theorem 4). In
particular, it follows that any arithmetic subgroup of a semi-simple 
algebraic group over ${\mathbb Q}$ has finite covolume in its real points.

The second chapter deals with several algebraic properties of arithmetic 
groups. As a consequence of reduction theory, we show that these groups are 
finitely generated. In fact they admit a finite presentation (Theorem 6). We 
give some explicit presentations of ${\rm SL}_N ({\mathbb Z})$, $N \geq 2$, 
and of the Chevalley groups $G({\mathbb Z})$ (5.6--5.8). We then show that, 
up to conjugation, arithmetic groups contain only finitely many finite 
subgroups (Theorem 7). Furthermore, they always contain a torsion free 
subgroup of finite index (Theorem 8). Following Minkowski, one can compute 
the least common multiple of the order of the finite subgroups of ${\rm GL}_N 
({\mathbb Z})$ (6.3). Coming back to $N=2$, we prove that any torsion free 
subgroup of ${\rm SL}_2 ({\mathbb Z})$ is a free group (Theorem 10). We 
conclude this section with the open problem, raised by Nahm, of finding the 
minimal index of torsion free subgroups of ${\rm SL}_N ({\mathbb Z})$, $N 
\geq 3$.

One of the main properties of arithmetic groups is their ``rigidity'' inside
the corresponding algebraic and Lie groups, at least when their rank is 
bigger than one. A lot of work has been accomplished on this theme. We start 
Chapter 3 with the congruence subgroup property, which states that any 
subgroup of finite index in $\Gamma$ contains the group of elements 
congruent to the identity modulo some integer. This property holds for 
arithmetic subgroups of simple simply connected Chevalley groups of rank 
bigger than one (Theorem 13), but it is wrong for ${\rm SL}_2$ (Corollary 
18). When studying that problem, Bass, Milnor and Serre discovered that, 
under suitable hypotheses, any linear representation of $\Gamma$ over 
${\mathbb Q}$ coincides with an algebraic representation on some subgroup of 
finite index (Proposition 14). This important rigidity property has many 
consequences, including the fact that the abelianization of $\Gamma$ is 
finite (Corollary 17).

Another approach to rigidity is Kazhdan's property (T), as explained in 
Theorems 19 and 20. Finally, we state the famous result of Margulis 
(Selberg's conjecture) that any discrete subgroup of 
finite covolume in a
simple, non-compact, connected Lie group of rank bigger than one is 
``arithmetic'' in a suitable sense (Theorem 21). This follows from a 
``superrigidity'' theorem for representations of arithmetic groups (Theorem 
22). Finally, we give another result of Margulis (Theorem 23), which states 
that arithmetic groups have rather few normal subgroups.

There are many results on arithmetic groups which are not covered by these 
notes. These include the different methods to compactify the quotient of a 
symmetric domain by the action of an arithmetic group (Baily-Borel-Satake, 
Borel-Serre$\ldots$), the cohomology of arithmetic groups (Borel, Serre, 
Franke,$\ldots$), and the ergodicity of their action on Lie groups 
(Margulis, Ratner,$\ldots$).

\medskip

I thank P. Cartier, B. Julia, W. Nahm, N. Nekrasov and J-P. Serre for 
helpful discussions.

\newpage

\centerline{\Large \bf I. Reduction theory}

\section{The reduction theory of quadratic forms}

\subsection{ \ } 

Groups of matrices with integral coefficients first appeared, in the work of 
Gauss, Hermite, Minkowski and others, as the  symmetry groups
for a
specific diophantine problem: which integers can be represented by a given 
quadratic form?

Recall that any positive integer is the sum of four squares. More generally, 
consider a quadratic form in $N$ variables
\begin{equation}
\label{soule:eq1}
\varphi (x) = \sum_{1 \leq i,j \leq N} a_{ij} \, x_i \, x_j
\end{equation}
where $a_{ij} = a_{ji}$. We assume that $\varphi$ is positive definite: for
any vector $x = (x_i) \in {\mathbb R}^N$, $\varphi (x) \geq 0$ and $\varphi 
(x) = 0$ iff $x=0$. When all coefficients $a_{ij}$ are integers, we say that 
a given integer $k \in {\mathbb N}$ is {\it represented by} $\varphi$ if
there exists $x \in {\mathbb Z}^N$ such that $\varphi (x) = k$.

Let now $\gamma \in {\rm GL}_N ({\mathbb Z})$ be an $N$ by $N$ square matrix 
with integral coefficients, the inverse of which has integral coefficients 
as well, i.e. such that $\det (\gamma) = \pm 1$. Let $^t{\gamma}$ be
the transpose of the matrix $\gamma$. When we change the coordinates of $x$ 
by $^t{\gamma}$, we get a new quadratic form $\gamma \cdot \varphi$:
\begin{equation}
\label{soule:eq2}
(\gamma \cdot \varphi) (x) = \varphi \, (^t{\gamma} (x)) \, ,
\end{equation}
for all $x \in {\mathbb R}^N$. Since $\gamma_1 (\gamma_2 (\varphi)) = 
(\gamma_1 \, \gamma_2) (\varphi)$, formula (\ref{soule:eq2}) defines an action of 
${\rm GL}_N ({\mathbb Z})$ on positive definite quadratic forms. It follows
from (\ref{soule:eq2}) that $k$ is represented by $\varphi$ iff $k$ is represented 
by $\gamma \cdot \varphi$ for all $\gamma \in {\rm GL}_N ({\mathbb Z})$. 
Therefore, when studying the integral values of $\varphi$ we may 
replace $\varphi$
by any form equivalent to it.

The {\it reduction theory of quadratic forms} consists in studying the 
orbits of ${\rm GL}_N ({\mathbb Z})$ on quadratic forms, and finding a good 
set of representatives for this action. Let $X$ be the set of positive 
definite quadratic forms 
$$
\varphi (x) = \sum_{1 \leq i,j \leq N} a_{ij} \, x_i \, x_j \,
,
$$
where $a_{ij} = a_{ji}$ are real numbers. We look for a ``small'' subset $D 
\subset X$ such that any point of $X$ is the translate of a point in $D$
by an element of $\Gamma$; in other words:
$$
\Gamma \cdot D = X \, .
$$

\subsection { \ }

Consider first the case $N=2$. Any positive definite binary form $\varphi$ 
can be written uniquely
\begin{equation}
\label{soule:eq3}
\varphi (x,y) = a(zx + y) (\bar z x + y)
\end{equation}
where $a > 0$ and $z$ is a complex number with positive imaginary part. The 
action of positive scalars on $X$ commutes with ${\rm GL}_N ({\mathbb Z})$ 
(for any $N \geq 2$) and (\ref{soule:eq3}) tells us that, when $N=2$, the quotient 
$X / {\mathbb R}_+^*$ is isomorphic to the Poincar\'e upper half plane
$$
{\mathcal H} = \{ z \in {\mathbb C} \vert \, {\rm Im} (z) > 0 \} \, .
$$
The action of $\gamma = \begin{pmatrix} a& b \\ c &d \end{pmatrix} \in {\rm 
SL}_2 ({\mathbb Z})$ is given by
\begin{equation}
\label{soule:eq4}
\gamma (z) = \frac{az + b}{cz + d} \, .
\end{equation}

\medskip

\noindent {\bf Theorem 1.} {\it Let $D$ be the set of $z \in {\mathcal H}$ 
such that $\vert z \vert \geq 1$ and $\vert {\rm Re} (z) \vert \leq 1/2$ 
(Figure 1). Then}
$$
{\mathcal H} = \Gamma \cdot D \, .
$$

\medskip

\noindent {\bf Remark.} If $z$ lies in the interior of $D$ and $\gamma (z) = 
z$, then $\gamma = \pm \, {\rm Id}$. Furthermore, if $\vert z \vert > 1$ and 
$\gamma (z) = z$, then ${\rm Re}  (z) = \pm \, 1/2$ and $\gamma (z) = 
 z+1$ or $z-1$.

\medskip

\noindent {\it Proof} (see \cite{soule:Se1}, VII, 1.2). Fix $z \in {\mathcal H}$.
We have
$$
{\rm Im} (\gamma (z)) = \frac{{\rm Im} (z)}{\vert cz + d \vert^2} \, ,
$$
and, given $A > 0$, there exist only finitely many $(c,d) \in {\mathbb Z}^2$ 
such that $\vert cz + d \vert^2 \leq A$. Therefore we can choose $\gamma$ 
such that ${\rm Im} (\gamma (z))$ is maximal. Let $T = \begin{pmatrix} 1 &1 
\\ 0 &1 \end{pmatrix}$. Since $T(z) = z+1$ we can choose $n \in {\mathbb Z}$ 
such that
$$
\vert {\rm Re} \, (T^n \, \gamma (z)) \vert \leq 1/2 \, .
$$
We claim that $z' = T^n \, \gamma (z)$ lies in $D$, i.e. $\vert z' 
\vert \geq 1$. Indeed, if $S = \begin{pmatrix} 0 &-1 \\ 1 &0 \end{pmatrix}$, 
we get $S(z') = -1/z'$, hence
$$
{\rm Im} (S(z')) = \frac{{\rm Im} (z')}{\vert z' \vert^2} \, .
$$
Since the imaginary part of $z'$ is maximal, this implies $\vert z' \vert 
\geq 1$.

$$
\includegraphics{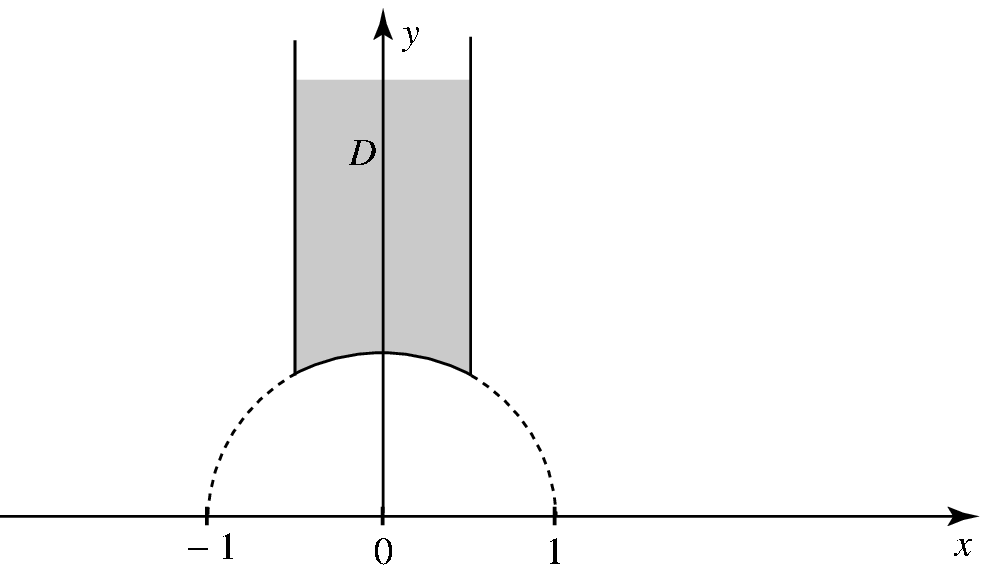}
\hbox{Figure 1}
$$

\section{Siegel sets}

More generally, when $N \geq 2$, any positive definite real quadratic form 
$\varphi$ in $N$ variables can be written uniquely in the following way:
\begin{eqnarray}
\label{soule:eq5}
\varphi (x) &= &t_1 (x_1 + u_{12} \, x_2 + u_{13} \, x_3 + \cdots + u_{1N} 
\, x_N)^2 \\
&+ &t_2 (x_2 + u_{23} \, x_3 + \cdots + u_{2N} \, x_N)^2 \nonumber \\
&+ & \cdots \nonumber \\
&+ &t_N \, x_N^2 \, , \nonumber
\end{eqnarray}
where $t_1 , \ldots , t_N$ are positive real numbers and $u_{ij} \in 
{\mathbb R}$.

\medskip

\noindent {\bf Theorem 2} (\cite{soule:B}, Th. 1.4). {\it After replacing 
$\varphi$ by $\gamma \cdot \varphi$ for some $\gamma \in {\rm GL}_N 
({\mathbb Z})$, we can assume that
$$
\vert u_{ij} \vert \leq 1/2 \quad \hbox{when} \quad 1 \leq i < j < N
$$
and
$$
t_i \leq \frac{4}{3} \, t_{i+1} \quad \hbox{when} \quad 1 \leq i \leq N-1 \,
.
$$
}

\medskip

The subset ${\mathfrak S} $ of $X$ defined by the inequalities of Theorem 2 is
called a {\it Siegel set} (see (\ref{soule:eq8}) below for a general definition).
When $N=2$, ${\mathfrak S} $ is the shaded region in Figure 2 below, and
Theorem 2 follows from Theorem 1.
$$
\includegraphics{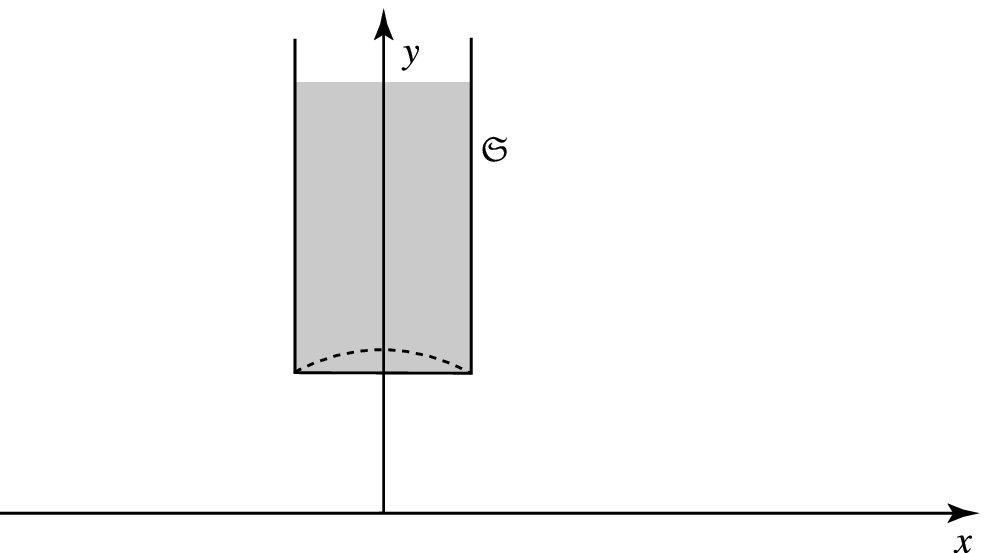}
\hbox{Figure 2}
$$

\section{Arithmetic groups}

\subsection{ \ }

Let $N \geq 1$ be an integer and $G$ a subgroup of ${\rm GL}_N ({\mathbb 
C})$. The group $G$ is called {\it linear algebraic over} ${\mathbb Q}$ if 
there exist polynomials $P_1 , \ldots , P_k$ with coefficients in ${\mathbb 
Q}$ in the variables $x_{ij}$, $1 \leq i,j \leq N$ and $u$ such that $G$ is 
the set of elements $g = (g_{ij}) \in {\rm GL}_N ({\mathbb C})$ such that
$$
P_1 (g_{ij} , \det (g)^{-1}) = P_2 (g_{ij} , \det (g)^{-1}) = \cdots = P_k 
(g_{ij} , \det (g)^{-1}) = 0 \, .
$$
The group $G({\mathbb Q}) = G \cap {\rm GL}_N ({\mathbb Q})$ is called the 
group of {\it rational points} of $G$.

Given $\Gamma_1$ and $\Gamma_2$ two subgroups of $G$, we say that $\Gamma_1$ 
and $\Gamma_2$ are {\it commensurable} when their intersection $\Gamma_1 
\cap \Gamma_2$ has finite index in both $\Gamma_1$ and $\Gamma_2$.

\medskip

\noindent{\bf Definition.} {\it Given $N \geq 1$ and $G \subset {\rm GL}_N 
({\mathbb C})$ a linear algebraic group over ${\mathbb Q}$, an {\rm 
arithmetic subgroup} of $G$ is a subgroup $\Gamma$ of $G({\mathbb Q})$ which 
is commensurable with $G \cap {\rm GL}_N ({\mathbb Z})$.}

\subsection{ \ }

A morphism $f : G \to G'$ of linear algebraic groups over ${\mathbb Q}$ is a 
group morphism defined by polynomials with coefficients in ${\mathbb Q}$ 
(note that $f$ needs {\it not} extend to a morphism between the ambient linear
groups).

\medskip

\noindent{\bf Proposition 3} (\cite{soule:B} Cor. 7.13, (3)). {\it If $\Gamma 
\subset G({\mathbb Q})$ is an arithmetic subgroup of $G$ and $f : G \to G'$ 
a morphism of linear algebraic groups over ${\mathbb Q}$, the image 
$f(\Gamma)$ is contained in some arithmetic group $\Gamma' \subset G' 
({\mathbb Q})$.}

\medskip

\noindent{\bf Remarks.} 1) If $G \subset {\rm GL}_N ({\mathbb C})$ is a 
linear algebraic group over ${\mathbb Q}$, we may consider its ring of 
rational functions
$$
A = {\mathbb Q} \, [x_{ij} , u] / \langle P_1 , \ldots , P_k \rangle \, .
$$
This ${\mathbb Q}$-algebra is finitely generated over ${\mathbb Q}$ and 
carries a Hopf structure coming from the group structure of $G$. Therefore 
$G_{{\mathbb Q}} = {\rm Spec} (A)$ is an {\it affine group scheme} over 
${\mathbb Q}$ [W]. The group $G$ is the set of complex points $G = {\rm Hom}
({\rm Spec} ({\mathbb C}) ,$ $G_{{\mathbb Q}})$ and $G({\mathbb Q})$ is the
set of rational points of $G_{{\mathbb Q}}$. Note that the definition of 
$G_{{\mathbb Q}}$ does not refer anymore (up to isomorphism) to a particular 
linear embedding.

\smallskip

\noindent 2) When $f : G \to G'$ is an isomorphism, it follows from 
Proposition 3 that $f(\Gamma)$ is arithmetic. This proves that the class of 
arithmetic subgroups of $G$ is intrinsic, i.e. it depends only on 
$G_{{\mathbb Q}}$ and not on the choice of the embedding $G \subset {\rm 
GL}_N ({\mathbb C})$.

\smallskip

\noindent 3) Another consequence of Proposition 3 is that, for any lattice 
$\Lambda \subset {\mathbb Q}^N$ (i.e. a free ${\mathbb Z}$-module of 
rank $N$), the group $\Gamma$ of elements $\gamma \in G ({\mathbb Q})$ such 
that $\gamma (\Lambda) = \Lambda$ is arithmetic.

\subsection{ \ }

The following general construction of arithmetic groups is due to Chevalley.

 \bigskip

\noindent{\bf 3.3.1}
 
 \bigskip

Let $E$ be a finite dimensional real euclidean vector space, and $\Phi
\subset E$ a {\it root system} (\cite{soule:H}, 9.2). Let $L_0 \subset E$ be the 
lattice spanned by $\Phi$ (the {\it lattice of roots}) and $L_1$ the {\it
lattice of weights}, i.e. those $\lambda \in E$ such that $\langle 
\lambda , \alpha \rangle \in {\mathbb Z}$ for all $\alpha \in \Phi$. The 
lattice $L_0$ is contained in $L_1$. Choose a lattice $L$ such that
$$
L_0 \subset L \subset L_1 \, .
$$

Given $\Phi$ and $L$, Chevalley defines as follows a linear algebraic group
$G$ over ${\mathbb Q}$ (\cite{soule:C}, \cite{soule:St}, \cite{soule:H} Chapter VII). Let
${\mathcal L}$ be a complex Lie algebra and ${\mathcal H} \subset {\mathcal L}$
a Cartan subalgebra such that $\Phi$ is the set of roots of ${\mathcal L}$ 
with respect to ${\mathcal H}$. If $\ell = \dim_{{\mathbb C}} {\mathcal H}$ 
and if $\Delta = \{ \alpha_1 , \ldots , \alpha_{\ell} \}$ is a basis of
$\Phi$, we can choose a {\it Chevalley basis} of ${\mathcal L}$, i.e. 
a basis $\{ X_{\alpha} , \alpha \in \Phi \, ; \ H_i , 1 \leq i \leq \ell \}$ 
such that $H_i \in {\mathcal H}$ and
$$
[H_i , H_j] = 0
$$
$$
[H_i , X_{\alpha}] = \langle \alpha , \alpha_i \rangle \, X_{\alpha}
$$
$$
[X_{\alpha} , X_{-\alpha}] \in \bigoplus_{i=1}^{\ell} {\mathbb Z} \, H_i
$$
$$
[X_{\alpha} , X_{\beta}] = \left\{ \begin{matrix} N_{\alpha\beta} \, 
X_{\alpha + \beta} &\hbox{when} \ \alpha + \beta \in \Phi \hfill \\ 0
\hfill &\hbox{otherwise,} \ \alpha + \beta \ne 0 \, . \end{matrix} \right. 
$$
Here $N_{\alpha\beta} \in {\mathbb Z}$ and $N_{-\alpha , -\beta} = 
-N_{\alpha\beta}$. Consider a faithful (i.e. injective) representation
$$
\rho : {\mathcal L} \to {\rm End} (V)
$$
of ${\mathcal L}$ on a finite dimensional
complex vector space $V$ such that $L$ is the set of
weights of $\rho$ (\cite{soule:H}, Ex. 21.5). 
For any root $\alpha \in \Phi$, the endomorphism
$\rho (X_{\alpha})^n = 0$ when $n$ is big enough so it makes sense to define 
$G$ as the group of endomorphisms of $V$ generated by the exponentials
$$
\exp (t \, \rho (X_{\alpha})) = \sum_{n \geq 0} t^n \, \frac{\rho 
(X_{\alpha})^n}{n!}
$$
for all $t \in {\mathbb C}$ and $\alpha \in \Phi$.

One can choose a basis of $V$ such that its ${\mathbb Z}$-span $M$ is stable 
by the action of every endomorphism $\frac{\rho (X_{\alpha})^n}{n!}$, $n 
\geq 1$, $\alpha \in \Phi$ (such a lattice is called {\it admissible}). If
the embedding $G \subset {\rm GL}_r ({\mathbb C})$ is defined by 
such a basis ($r
= \dim_{{\mathbb C}} V$), $G$ is the set of zeroes of polynomials with 
${\mathbb Q}$-coefficients (\cite{soule:St}, \S~5, Th. 6). It can be shown 
\cite{soule:C} \cite{soule:St} that, up to canonical isomorphism, the linear algebraic 
group $G$ over ${\mathbb Q}$ depends only on $\Phi$ and $L$.

When $L = L_0$, the group $G$ is called {\it adjoint}, and when $L=L_1$ it is 
called {\it simply connected} (or ``universal'').

\bigskip

\noindent{\bf 3.3.2}

\bigskip

Let $\Phi$, $L$, $\rho$ and $M$ be as above. The group $G({\mathbb Z}) = \{ 
g \in G$ such that $g(M) = M\}$ is an arithmetic subgroup of $G$. Up to 
canonical isomorphism (defined by means of polynomials with integral 
coefficients, respecting the inclusion $G({\mathbb Z}) \subset G({\mathbb 
Q})$) it depends only on $\Phi$ and $L$. In fact, Chevalley proves in 
\cite{soule:C} that $(\Phi , L)$ defines an affine group scheme $G_{{\mathbb Z}}$ 
over ${\mathbb Z}$, and 
$$
G({\mathbb Z}) = {\rm Hom} ({\rm Spec} ({\mathbb Z}) , G_{{\mathbb Z}})
$$
is its set of integral points.

\subsection{ \ }

The group of integral points of the simply connected Chevalley 
group scheme of type $A_n$ (resp. $B_n$) is the group
 ${\rm SL}_n ({\mathbb Z})$ of integral matrices with determinant one
(resp. the group ${\rm Sp}_{2n} ({\mathbb Z})$ of symplectic matrices 
in ${\rm SL}_{2n} ({\mathbb
Z})$).

Another example is the simply connected Chevalley group scheme $G$ of type $E_7$ 
over ${\mathbb Z}$ and its set $G({\mathbb Z})$ of integral points. Consider the 
split Lie group $E_{7(+7)}$ of type $E_7$ and its fundamental representation 
$E_{7(+7)} \subset {\rm Sp}_{56} ({\mathbb R})$ of dimension 56, as described
 in \cite{soule:CJ}, Appendix B. Let
  $E_7 ({\mathbb Z}) = E_{7 (+7)} \cap {\rm Sp}_{56} ({\mathbb Z})$ as in \cite{soule:HT}.
  We shall prove in the Appendix that
\begin{equation}
E_7 ({\mathbb Z}) = G({\mathbb Z}) \, .
\end{equation}

\subsection{ \ } 

Assume $F$ is a number field (a finite extension of ${\mathbb Q}$). We can 
define linear algebraic groups $G$ over $F$ in the same way as when $F = 
{\mathbb Q}$, by choosing a complex embedding of $F$ or more intrinsically 
as in 3.2 Remark 1). Matrices in $G$ with coefficients in the integers of 
$F$ are arithmetic groups.

However, these definitions do {\it not} enlarge the class of arithmetic 
groups. Indeed $G(F)$ can be viewed as the group of rational points 
$H({\mathbb Q}) = G(F)$, where $H = {\rm Res}_{F/{\mathbb Q}} \, G$ is the {\it 
restriction of scalars} of $G$ from $F$ to ${\mathbb Q}$ (\cite{soule:B} 7.16, 
\cite{soule:OV} Chapter 3, \S~6.1).

For example, let $d > 0$ be a positive integer which is not a square. 
Consider the subgroup $H$ of ${\rm GL}_2 ({\mathbb C})$ made of matrices $g 
= (g_{ij})$ such that $g_{11} = g_{22}$, $g_{21} = dg_{12}$ and $\det (g) = 
1$. In other words each $g \in H$ can be written
$$
g = x \cdot 1 + y \cdot \sigma
$$
where $\sigma = \begin{pmatrix} 0 &1 \\ d &0 \end{pmatrix}$. Note that 
$\sigma^2 = d \cdot 1$. The group $H$ is the restriction of scalars ${\rm 
Res}_{F/{\mathbb Q}} \, {\rm GL}_1$ where $F = {\mathbb Q} \, (\sqrt d)$. Note 
that $H({\mathbb R})$ is isomorphic to ${\mathbb R}^*$ (map $x \cdot 1 + y 
\cdot \sigma$ to $x+y \, \sqrt d$) but $H$ is not isomorphic to ${\rm GL}_1$
over ${\mathbb Q}$. We have $H({\mathbb Q}) = F^*$, and the group of units 
${\mathcal O}_F^*$ is an arithmetic subgroup of $H$.

\section{The reduction theory of arithmetic groups}

\subsection{ \ }

Theorem 1 can be extended to all arithmetic subgroups of reductive groups.
We first need some definitions. A linear algebraic group $U$ is called {\it
unipotent} (resp. {\it solvable}) when there exists a finite filtration 
$ \cdots \subset U_i
\subset U_{i+1} \subset \cdots \subset U$ of $U$ by (Zariski) closed normal
subgroups such that each quotient $U_{i+1} / U_i$ is isomorphic (over 
${\mathbb Q}$) to the additive group (resp. is abelian). If $G$ is any 
linear algebraic group over ${\mathbb Q}$, the {\it unipotent radical} $R_u 
(G)$ (resp. the {\it radical} $R(G)$) is the maximal closed connected 
unipotent (resp. solvable) normal subgroup of $G$. Of course $R_u (G)$ is 
contained in $R(G)$. The group $G$ is called {\it reductive} (resp. {\it 
semi-simple}) when $R_u (G) = \{ 1 \}$ (resp. $R(G) = \{ 1 \}$).

Let $G$ be a reductive linear algebraic group over ${\mathbb Q}$, let $G^0$ 
be the connected component of the unit element $1 \in G$, and let $P$ be a 
minimal {\it parabolic subgroup} of $G^0$ over ${\mathbb Q}$, i.e. a 
minimal closed connected subgroup $P \subset G^0$ such that the variety 
$G/P^0$ is projective. According to \cite{soule:B}, Th. 11.4, i), one can write 
$P$ as a product of subgroups
\begin{equation}
\label{soule:eq6}
P = M \cdot S \cdot U \, ,
\end{equation}
where $U = R_u (P)$ and $S$ is a maximal split ${\mathbb Q}$-torus of $P$ 
(i.e. $S$ is isomorphic over ${\mathbb Q}$ to a power of the 
multiplicative group). Let $X(Z(S))$ be the set of characters $\chi : Z(S) 
\to {\rm GL}_1$ over ${\mathbb Q}$ of the centralizer $Z(S)$ of $S$ in 
$G^0$. Then $M$ is defined as the connected component of $1$ in the 
intersection $\cap \ker (\chi)$ of the kernels of
all the characters $\chi \in X (Z(S))$.

Furthermore, if $G({\mathbb R}) = G \cap {\rm GL}_N ({\mathbb R})$ is the 
group of real points of $G$ and $K$ a maximal compact subgroup of this Lie 
group $G({\mathbb R})$, we may write, according to \cite{soule:B} 11.19,
\begin{equation}
\label{soule:eq7}
G({\mathbb R}) = P({\mathbb R}) \cdot K = M^0 \cdot N \cdot A \cdot K
\end{equation}
where $M^0$ (resp. $A$) is the (usual) connected component of $1$ in 
$M({\mathbb R})$ (resp. $S({\mathbb R})$), and $N = U({\mathbb R})$. The 
decomposition (\ref{soule:eq7}) generalizes the Iwasawa decomposition.

\subsection{ \ }

For example, when $G = {\rm GL}_N ({\mathbb C})$, we can choose for $K$ the 
group $O_N ({\mathbb R})$ of orthogonal matrices, for $P$ lower triangular
matrices, for $S$ diagonal ones,  for
$N$ lower unipotent matrices and $M = \{ 1 \}$.
Define a map from ${\rm GL}_N ({\mathbb R})$ to the space $X$ of real 
positive definite quadratic forms by the formula
$$
\varphi (x) = \Vert ^tg (x) \Vert^2 \, .
$$
Using this map, we see that (\ref{soule:eq7}) follows from (\ref{soule:eq6}) in this 
case.

\subsection{ \ }

We come back to the notations of \S~4.1. Let $X(S) = {\rm Hom}_S (S,{\rm 
GL}_1)$ be the set of characters of $S$ over ${\mathbb Q}$, and let $\Phi
\subset X(S)$ be the set of roots of $G$. The group $U$ defines an ordering of $\Phi$
(\cite{soule:B}, 11.6 (3)) and we let $\Delta \subset \Phi$ be the set of positive 
simple roots. For any real number $t > 0$ let
$$
A_t = \{ a \in A \vert \, \alpha (a) \leq t \ \hbox{for all} \ \alpha \in \Delta \}
\, .
$$
If $\omega$ is any compact neighbourhood of $1$ in $M^0 \cdot N$ we define
\begin{equation}
\label{soule:eq8}
{\mathfrak S} _{t,\omega} = \omega \cdot A_t \cdot K \, ,
\end{equation}
a subset of $G({\mathbb R})$ by (\ref{soule:eq7}). This set
 ${\mathfrak S}_{t,\omega}$ is called a {\it Siegel set}.

\medskip

\noindent {\bf Theorem 4.} {\it Let $G$ be a reductive linear algebraic 
group over ${\mathbb Q}$ and $\Gamma$ an arithmetic subgroup of $G$.}
\begin{itemize}
\item[i)] (\cite{soule:B}, Th. 1.3.1) {\it One can choose $t > 0$ and $\omega$ as 
above, and a finite subset $C$ in $G({\mathbb Q})$ such that}
$$
G({\mathbb R}) = \Gamma \cdot C \cdot {\mathfrak S} _{t,\omega} \, .
$$
\item[ii)] (\cite{soule:B}, Th. 15.4) {\it For any choice of $t$ and $\omega$, and 
any $g \in G({\mathbb Q})$, the set of elements $\gamma \in \Gamma$ such 
that $g \cdot {\mathfrak S} _{t,\omega}$ meets
$\gamma \cdot {\mathfrak S}{t,\omega}$ is} finite.
\item[iii)] (\cite{soule:B}, Lemma 12.5) {\it The Haar measure of any Siegel set 
is finite iff the set of rational characters $X(G^0)$ is trivial.}
\end{itemize}

\medskip

\noindent {\bf Corollary 5.} {\it With the same hypotheses:}
\begin{itemize}
\item[i)] {\it The quotient $\Gamma \backslash G({\mathbb R})$ is compact 
iff $G$ is} anisotropic {\it over ${\mathbb Q}$ {\rm (i.e.} $S=U=\{ 1 \})$.}
\item[ii)] {\it The (invariant) volume of $\Gamma \backslash G({\mathbb R})$ 
is finite iff $G^0$ has no nontrivial cha\-racter over ${\mathbb Q}$ (
e.g. if $G$ is semi-simple).}
\end{itemize}

\newpage

\centerline{\Large \bf II. Some algebraic properties}

\centerline{\Large \bf of arithmetic groups}

\section{Presentations}

\subsection{ \ }

Let $S$ be a set. The {\it free group} $F(S)$ over $S$ is defined by the 
following universal property: $S$ is contained in $F(S)$ and, given any 
group $G$ and any map of sets $\varphi : S \to G$, there exists a unique 
group morphism $\tilde\varphi : F(S) \to G$ which coincides with $\varphi$ on
$S$, i.e. such that the diagram
$$
\begin{matrix}
F(S) &\overset{\tilde\varphi}{\longrightarrow} &G \\
\hfill \nwarrow &&\!\!\!\!\!\!\nearrow_{\varphi} \\
&\!\!\!\!\!\!S
\end{matrix}
$$
commutes.

Clearly $F(S)$ is unique up to unique isomorphism. A construction of $F(S)$ 
is given in \cite{soule:L} I, \S~8, Prop. 7.

\subsection{ \ }

Given a group $G$, a {\it presentation} of $G$ is a pair $(S,R)$ where $S 
\subset G$ is a subset of $G$ and $R \subset F(S)$ is a subset of the free 
group over $S$ such that
\begin{itemize}
\item[i)] $S$ spans $G$, i.e. the canonical map $F(S) \to G$ is 
surjective;
\item[ii)] the kernel of the map $F(S) \to G$ is the smallest normal 
subgroup $\langle R \rangle$ of $F(S)$ containing $R$.
\end{itemize}

\medskip

It follows from i) and ii) that $G$ is isomorphic to $F(S) / \langle R 
\rangle$. We say that $G$ is generated by $S$, with relations $r=1$ for all 
$r \in R$.

When $S$ and $R$ are finite, $(S,R)$ is a {\it finite presentation} of $G$.

\subsection{ \ }

For example, when $S = \{ x,y \}$ consists of two elements and $R = \{ x^2 , 
y^2 , (xy)^3 \}$ the group $G = F(S) / \langle R \rangle$ is the group $S_3$ 
of permutations of three elements. This can be seen by 
mapping $x$ (resp. $y$)
to the permutation $(123) \to (213)$ (resp. $(123) \to (132)$).

\subsection{ \ }

Here are two finite presentations of ${\rm SL}_2 ({\mathbb Z})$:
\begin{itemize}
\item[a)] ${\rm SL}_2 ({\mathbb Z})$ is generated by 
$x = \begin{pmatrix} 1 &1
\\ 0 &1 \end{pmatrix}$ and $y = \begin{pmatrix} 1 &0 \\ 1 &1 \end{pmatrix}$ 
with relations
$$
(x^{-1} \, y \, x^{-1})^4 = 1 \qquad \hbox{and} \qquad x \, y^{-1} \, x = 
y^{-1} \, x \, y^{-1} \, .
$$
\item[b)] ${\rm SL}_2 ({\mathbb Z})$ is generated by $W = \begin{pmatrix} 0 
&-1 \\ 1 &0 \end{pmatrix}$ and $A = \begin{pmatrix} 1 &-1 \\ 1 &0 
\end{pmatrix}$ with relations $W^2 = A^3$ and $W^4 = 1$.
\end{itemize}

\subsection{ \ }

In general we have the following

\medskip

\noindent {\bf Theorem 6} (\cite{soule:B2}, \cite{soule:R}). {\it Let $\Gamma$ be an 
arithmetic subgroup of a linear algebraic group $G$ over ${\mathbb Q}$. 
Then $\Gamma$ is finitely presented.}

\medskip

In other words, $\Gamma$ admits a finite presentation. Let us indicate why
$\Gamma$ is finitely generated. Let $K$ be a maximal compact subgroup of 
$G({\mathbb R})$ and $X = G({\mathbb R}) / K$. We claim that we can find a 
closed subset $D \subset X$ such that
\begin{itemize}
\item[a)] $\Gamma \cdot D = X$; 
\item[b)] the subset $S\subset \Gamma$  of those $\gamma$ such that
$\gamma\cdot D \cap D \ne \emptyset$ is finite.
\end{itemize}

Indeed, when $G$ is reductive we can take for $D$ the union $\underset{g \in 
C}{\bigcup} g \cdot {\mathfrak S} _{t,\omega}$ of finitely many translates of 
a (well chosen) Siegel set; see Theorem 4, i) and ii). When $G$ is 
arbitrary, it is a semi-direct product
\begin{equation}
\label{soule:eq9}
G = R_u (G) \cdot H,
\end{equation}
where $H$ is reductive over ${\mathbb Q}$ (\cite{soule:B} 7.15) and $R_u (G)$ is 
the unipotent radical of $G$. The quotient $(R_u (G) \cap \Gamma) \backslash 
R_u (G)$ is compact so we can choose a compact subset $\Omega \subset R_u 
(G)$ such that $R_u (G) = (R_u (G) \cap \Gamma) \cdot \Omega$. Let $D' 
\subset H({\mathbb R}) / K$ be such that a) and b) are true for $D'$ and
$\Gamma \cap H$ (note that $K \cap R_u (G)$ is trivial). Then one can check 
that $D = \Omega \cdot D'$ and $\Gamma$ satisfy a) and b).

From this we derive that $S$ spans $\Gamma$. Indeed, it follows from 
(\ref{soule:eq7}) and (\ref{soule:eq9}) that $X$ is homeomorphic to a euclidean space 
and, in particular, it is path connected. Given $\gamma \in \Gamma$, choose 
a continuous path $c : [0,1] \to X$ such that $x = c(0)$ lies in $D$ and 
$c(1) = \gamma \cdot x$. Since $c([0,1])$ is compact, there exists a finite 
sequence $\gamma_1 , \ldots , \gamma_k$ in $\Gamma$ such that $\gamma_1 = 
1$, $\gamma_i \cdot D \cap \gamma_{i+1} \cdot D \ne \emptyset$ when $i < k$,
and $\gamma_k = \gamma$. Define $s_i = \gamma_i^{-1} \, \gamma_{i+1}$, $i < 
k$. Since $s_i \cdot D \cap D \ne \emptyset$, these elements lie in $S$. On the
other hand,
$$
\gamma = s_1 \ldots s_{k-1} \, ,
$$
therefore $S$ spans $\Gamma$.

\subsection{ \ }

Theorem 6 gives us  general information, but it is still of interest to
have explicit presentations of arithmetic groups. For instance, consider the 
case of ${\rm SL}_N ({\mathbb Z})$, $N \geq 3$. For any pair of indices 
$(i,j)$, $1 \leq i \ne j \leq N$, denote by $x_{ij} \in {\rm SL}_N ({\mathbb 
Z})$ the matrix which is equal to one on the diagonal and at the entry 
$(i,j)$, and zero otherwise:
$$
x_{ij} =
\begin{pmatrix}
1&&0&&&j&& \\
&1&&0&&|&& \\
&&1&&&1& -i \\
&0&&\ddots &&0 \\
&&0&&1&\\
&&&0&&1\\
&&&&&&1
\end{pmatrix} \, .
$$
These matrices $x_{ij}$, $1 \leq i \ne j \leq N$, generate ${\rm SL}_N 
({\mathbb Z})$, and the following relations give a presentation:
$$
[x_{ij} , x_{k\ell}] = 1 \qquad \hbox{if} \ j \ne k \ \hbox{and} \ i \ne 
\ell \, ;
$$
$$
[x_{ij} , x_{jk}] = x_{ik} \qquad \hbox{if} \ i,j,k \ \hbox{are distinct} \, 
;
$$
$$
(x_{12} \, x_{21}^{-1} \, x_{12})^4 = 1 \, .
$$
As usual, $[g,h]$ is the commutator $ghg^{-1} h^{-1}$. This fact is due to 
Magnus and Nielsen (see \cite{soule:Mi} Cor. 10.3).

\medskip

\noindent {\bf Remark.} There are known bounds for the number of elementary
matrices $x^a_{ij}, a\in {\mathbb Z}$,  needed to write any
element of ${\rm SL}_N ({\mathbb
Z})$ \cite{soule:Ca}. For instance, any element of ${\rm SL}_3
({\mathbb Z})$  is the product of at most 60 elementary matrices.

\subsection{ \ }

The group ${\rm SL}_N ({\mathbb Z})$ can be generated by two elements only, 
for instance $x_{21}$ and the matrix $(g_{ij})$ where
$$
g_{ij} = \left\{ \begin{matrix}
1 \hfill &\hbox{if} \ 1 \leq i \leq N-1 \ \hbox{and} \ j=i+1 \\
(-1)^N &\hbox{if} \ (i,j) = (N,1) \hfill \\
0 \hfill &\hbox{otherwise} \hfill
\end{matrix}
\right.
$$
(\cite{soule:C-M}, p.83).
For a (long) list of defining relations between these two matrices, see
\cite{soule:C-M}, p.85.

\subsection{ \ }

Let $\Phi$ be a root system and $L$ a lattice such that $L_0 \subset L 
\subset L_1$ as in 3.3.1. Let $G({\mathbb Z})$ be the associated arithmetic 
group (3.3.2).

Choose a faithful representation $\rho : {\mathcal L} \to {\rm End} (V)$ 
with weight lattice $L$ as in 3.3.1. The group $G({\mathbb Z})$ is then 
generated by the endomorphisms
$$
x_{\alpha} = \exp (\rho(X_{\alpha})) \in {\rm End} (V)
$$
(\cite{soule:St}, Th. 18, Cor. 3, Example).

Assume furthermore that $\Phi$ is irreducible, $\Phi \ne A_1$, and $L = L_1$ 
(so that the Chevalley group $G$ is simple, simply connected and different 
from ${\rm SL}_2$). The following relations define $G({\mathbb Z})$ (and 
generalize 5.6) (\cite{soule:Be}, Satz 3.1):
$$
[x_{\alpha} , x_{\beta}] = \prod_{i,j} x_{i\alpha + j\beta}^{N(\alpha , 
\beta ; i,j)} \quad \hbox{when} \ \alpha + \beta \ne 0 \, ;
$$
$$
(x_{\alpha}^{-1} \, x_{-\alpha} \, x_{\alpha}^{-1})^4 = 1 \qquad \hbox{for 
any simple root $\alpha$.}
$$

Here $i$ and $j$ run over positive integers and the integers $N(\alpha , 
\beta ; i,j)$ are almost all zero ($N(\alpha , \beta ; 1,1) = 
N_{\alpha\beta}$ are the constants defining the Chevalley basis in 3.3.1).

\section{Finite subgroups}

\subsection{ \ }

\noindent {\bf Theorem 7} (\cite{soule:B2} \cite{soule:R}). 
{\it Let $\Gamma$ be an arithmetic
subgroup of a linear algebraic group $G$ over ${\mathbb Q}$. Up to 
conjugation, $\Gamma$ contains only finitely many finite subgroups.}

\medskip

\noindent {\it Proof.} Let $X = G({\mathbb R}) / K$ and $D \subset X$ be as 
in the proof of Theorem 6. Any finite subgroup $F \subset \Gamma$ is 
contained in a maximal compact subgroup $K'$ of $G({\mathbb R})$. Since $K' 
= g \, K \, g^{-1}$ is conjugate to $K$, the point $x = g \, K$ in $X$ is 
fixed by all $\gamma \in F$.

Let $y \in D$ and $\gamma' \in \Gamma$ be such that $x = \gamma' (y)$. Then, 
for all $\gamma \in F$, we have
$$
{\gamma'}^{-1} \, \gamma \, \gamma' (y) = y \, .
$$
In particular ${\gamma'}^{-1} \, \gamma \, \gamma' (D)$ meets $D$ and 
${\gamma'}^{-1} \, \gamma \, \gamma'$ lies in the finite set $S$ (Theorem 6, 
b)). This proves our assertion.

\subsection{ \ }

\noindent {\bf Theorem 8.} {\it If $\Gamma$ is an arithmetic subgroup of 
$G$, there exists a subgroup of finite index $\Gamma' \subset \Gamma$ which 
is torsion free.}

\medskip

\noindent {\it Proof.} By definition (3.1), $\Gamma$ is commensurable with $G
\cap {\rm GL}_N ({\mathbb Z})$ for some embedding of $G$ in ${\rm GL}_N 
({\mathbb C})$. So it is enough to prove Theorem 8 for ${\rm GL}_N ({\mathbb 
Z})$. 

\medskip

It follows from the  following lemma.

\medskip

\noindent {\bf Lemma 9.} {\it Let $p \geq 3$ be a prime integer and $\Gamma$ 
the set of elements $\gamma \in {\rm GL}_N ({\mathbb Z})$ which are 
congruent to the identity modulo $p$. Then $\Gamma$ is a torsion free 
subgroup of ${\rm GL}_N ({\mathbb Z})$.}

\medskip

\noindent {\it Proof.} Clearly $\Gamma$ is a subgroup of ${\rm GL}_N 
({\mathbb Z})$. If it was not torsion free, it would contain an element of 
prime order, say $\ell > 1$, so there would exist a square matrix $m \in M_N 
({\mathbb Z})$ not divisible by $p$ and some integer $\alpha \geq 1$ such 
that
\begin{equation}
\label{soule:eq10}
(1 + p^{\alpha} \, m)^{\ell} = 1 \, .
\end{equation}
From the binomial formula, we deduce from (\ref{soule:eq10}) that
\begin{equation}
\label{soule:eq11}
\ell \, p^{\alpha} \, m = - \sum_{i=2}^{\ell} \begin{pmatrix} \ell \\ i 
\end{pmatrix} p^{\alpha i} \, m^i \, .
\end{equation}
When $\ell \ne p$, the exact power of $p$ dividing $\ell \, p^{\alpha} \, m$ 
is $p^{\alpha}$. But the right hand side of (\ref{soule:eq11}) is divisible by 
$p^{2\alpha}$, so we get a contradiction.

When $\ell = p$, $p^{\alpha + 1}$ is the exact power of $p$ dividing the 
left hand side of (\ref{soule:eq11}). When $2 \leq i < p$, $p$ divides 
$\begin{pmatrix} p \\ i \end{pmatrix}$, therefore $p^{2\alpha + 1}$ divides 
$\begin{pmatrix} p \\ i \end{pmatrix} p^{\alpha i}$. Finally, since $p \geq 
3$, $p^{\alpha p}$ is also divisible by $p^{2\alpha + 1}$. Therefore 
$p^{2\alpha + 1}$ divides the right hand side of (\ref{soule:eq11}) and we get 
again a contradiction.

\subsection{ \ }

From Lemma 9, Minkowski got some information on the order of the finite
subgroups of ${\rm GL}_N ({\mathbb Z})$ (\cite{soule:M} 212-218, \cite{soule:Bo} \S~7, 
Exercises 5-8). Indeed, when $p \geq 3$, any finite subgroup $F \subset {\rm 
GL}_N ({\mathbb Z})$ maps injectively into the quotient group ${\rm GL}_N 
({\mathbb Z} / p)$, the order of which is
$$
a(N,p) = (p^N - 1) (p^N - p) \cdots (p^N - p^{N-1}) \, .
$$
If $\ell$ is an odd prime, and if the reduction of $p$ modulo $\ell^2$ is a 
generator of $({\mathbb Z} / \ell^2 \, {\mathbb Z})^*$, the power of $\ell$
dividing  $a(N,p)$ is exactly $\ell^{r(\ell , N)}$ with
$$
r(\ell , N) = \left[ \frac{N}{\ell - 1} \right] + \left[ \frac{N}{\ell (\ell 
- 1)} \right] + \left[ \frac{N}{\ell^2 (\ell - 1)} \right] + \cdots \, ,
$$
where $[x]$ denotes the integral part of the real number $x$. Conversely, it 
can be shown (loc. cit.) that ${\rm GL}_N ({\mathbb Z})$ contains a finite 
subgroup of order $\ell^{r(\ell , N)}$.

The same results are true when $\ell = 2$ and 
$$
r(2,N) = N + \left[ \frac{N}{2} \right] + \left[ \frac{N}{4} \right] + 
\cdots \, .
$$
If we denote by $m(N)$ the product over all primes $\ell$ of $\ell^{r(\ell , 
N)}$, we conclude that $m(N)$ is the least common multiple of the 
cardinality of the finite subgroups of ${\rm GL}_N ({\mathbb Z})$. For 
instance
$$
m(2) = 24 \, , \ m(3) = 48 \, , \ m(4) = 5760 , \ldots
$$

\subsection{ \ }

Let us come back to ${\rm SL}_2 ({\mathbb Z})$.

\medskip

\noindent {\bf Theorem 10.} {\it Let $\Gamma \subset {\rm SL}_2 ({\mathbb 
Z})$ be any torsion free subgroup. Then $\Gamma$ is a free group.}

\medskip

\noindent {\it Proof.}  Let ${\mathcal H}$ be the Poincar\'e upper 
half-plane. Recall from Theorem 1 that ${\rm SL}_2 ({\mathbb Z})$ acts upon 
${\mathcal H}$ with fundamental domain the set $D$ of those $z \in {\mathcal 
G}$ such that $\vert z \vert \geq 1$ and $\vert {\rm Re} (z) \vert \leq 
1/2$.

The stabilizer in ${\rm SL}_2 ({\mathbb Z})$ of any $z \in {\mathcal H}$ is 
finite. Indeed ${\mathcal H} = {\rm SL}_2 ({\mathbb R}) / {\rm SO}_2 
({\mathbb R})$ hence the stabilizer of $z$ is the intersection of the 
discrete group ${\rm SL}_2 ({\mathbb Z})$ with a conjugate of the compact 
group ${\rm SO}_2 ({\mathbb R})$.
Since $\Gamma$ is torsion free, it acts freely on ${\mathcal H}$ (it has no
fixed point).

Let $D_0 \subset D$ be the set of points $z \in {\mathcal H}$ such that 
$\vert z \vert = 1$ and $\vert {\rm Re} (z) \vert \leq 1/2$, and
$$
Y = {\rm SL}_2 ({\mathbb Z}) \cdot D_0
$$
the union of the translates of $D_0$ under ${\rm SL}_2 ({\mathbb Z})$:

$$
\includegraphics{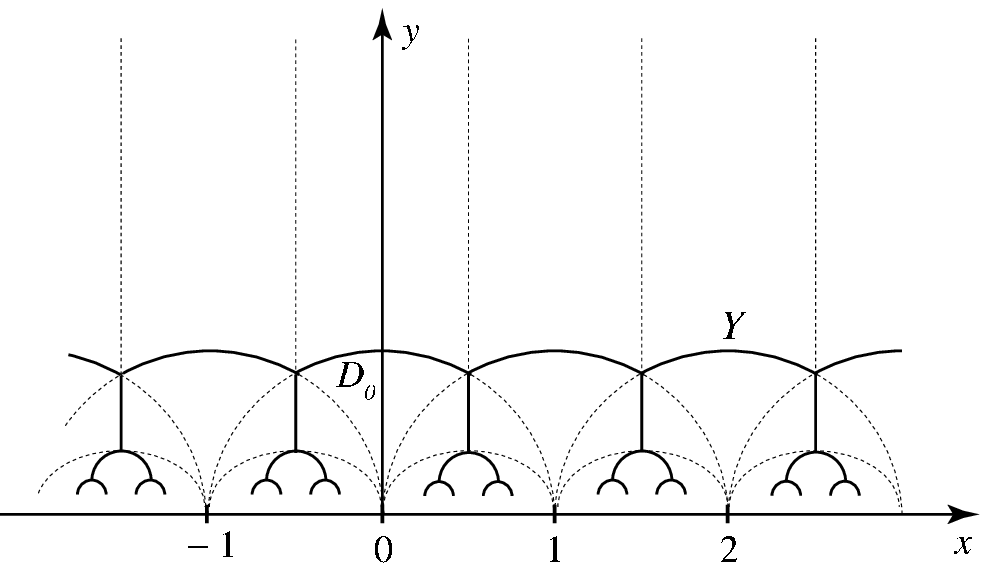}
\hbox{Figure 3}
$$

\medskip

\noindent {\bf Proposition 11} (\cite{soule:Se1} ). {\it The set $Y$ is (the
topological realization of) a tree.}

\medskip

\noindent {\it Proof of Proposition 11.} Clearly $Y$ is a graph, and we want 
to show that $Y$ can be contracted (deformed) to a point. Consider the 
retraction of $D$ onto $D_0$ which maps $z \in D$ to the point $z' \in D_0$ 
with the same abcissa as $z$. When $z \in D - D_0$ and $\gamma (z) \in D$ we
know that $\gamma (z) = z \pm 1$ (1.2, Remark). Therefore this retraction 
commutes with the action of ${\rm SL}_2 ({\mathbb Z})$ on $D$, and it can be 
extended to a retraction of ${\mathcal H} = {\rm SL}_2 ({\mathbb Z}) \cdot 
D$ onto $Y = {\rm SL}_2 ({\mathbb Z}) \cdot D_0$. Since ${\mathcal H}$ is
contractible to a point, the same is true for $Y$. q.e.d.

\bigskip

To end the proof of Theorem 10 note that $\Gamma$ acts freely on the tree 
$Y$, so it can be identified with the fundamental group of the quotient:
$$
\Gamma = \pi_1 (\Gamma \backslash Y) \, .
$$
This quotient $\Gamma \backslash Y$ is a connected graph and we have:

\medskip

\noindent {\bf Proposition 12.} {\it Let ${\mathcal X}$ be a connected 
graph. Then $\pi_1 ({\mathcal X})$ is free.}

\medskip

\noindent {\it Proof.} Choose a maximal tree $T \subset {\mathcal X}$. 
Clearly $T$ contains all the vertices of ${\mathcal X}$. Therefore, after 
contracting $T$, ${\mathcal X}$ becomes a ``bouquet'' of circles $B$. We get
$$
\pi_1 ({\mathcal X}) = \pi_1 (B) = F(S) \, ,
$$
where $S$ is the set of circles in $B$.

\subsection{ \ }

Let $\Gamma \subset {\rm SL}_2 ({\mathbb Z})$ be torsion free with finite 
index $e = [{\rm SL}_2 ({\mathbb Z}) : \Gamma]$. It can be shown that 12 
divides $e$ (see 6.6 below) and that the number of generators of $\Gamma$ is 
$1 + \frac{e}{12}$.

For instance, the subgroup of commutators
$$
\Gamma = [{\rm SL}_2 ({\mathbb Z}) , {\rm SL}_2 ({\mathbb Z})]
$$
has index 12 in ${\rm SL}_2 ({\mathbb Z})$. It is free on the two generators 
$\begin{bmatrix} 2 &1 \\ 1 &1 \end{bmatrix}$ and $\begin{bmatrix} 1 &1 \\ 1 
&2 \end{bmatrix}$.

\subsection{ \ }

Let $N \geq 2$ and let $\Gamma \subset {\rm SL}_N ({\mathbb Z})$ 
be a torsion free
normal subgroup of ${\rm SL}_N ({\mathbb Z})$. Any finite subgroup of ${\rm 
SL}_N ({\mathbb Z})$ maps injectively into the 
quotient ${\rm SL}_N ({\mathbb Z}) / \Gamma$. Therefore the index $[{\rm 
SL}_N ({\mathbb Z}) : \Gamma]$ is divisible by $m(N) / 2$, where $m(N)$ is 
as in 6.3. We have just seen that ${\rm SL}_2 ({\mathbb Z})$ contains a 
torsion free subgroup of index $m(2) / 2 = 12$. But, when $N \geq 3$, I do 
not know what  the minimal index of a torsion free subgroup of ${\rm
SL}_N ({\mathbb Z})$ is (a question raised by W.~Nahm).

\newpage

\centerline{\Large \bf III. Rigidity}

\section{The congruence subgroup problem}

\subsection{ \ }

Let $G \subset {\rm GL}_N ({\mathbb C})$ be a linear algebraic group over 
${\mathbb Q}$ and $\Gamma \subset G({\mathbb Q})$ an arithmetic subgroup. 
For any $a \geq  1$ we define a {\it congruence subgroup} $\Gamma (a)$ of
$\Gamma$. It consists of those matrices in $\Gamma \cap {\rm GL}_N ({\mathbb 
Z})$ which are congruent to the identity modulo $a$. This is a subgroup of 
finite index in $\Gamma$.

\medskip

\noindent {\bf Definition.} {\it We say that $G$ has property (CS) if any
$\Gamma' \subset \Gamma$ of finite index contains a congruence subgroup 
$\Gamma (a)$.}

\medskip

It can be shown that this property depends on $G$ only and neither on the
choice of $\Gamma$ nor on the embedding $G \subset {\rm GL}_N ({\mathbb
C})$.

\medskip

\noindent {\bf Theorem 13.} 
\begin{itemize}
\item[i)] (\cite{soule:BMS} \cite{soule:V} \cite{soule:Ma}) {\it If $\Phi$ is an irreducible 
root lattice different from $A_1$ and if $L = L_1$, the (simple and simply 
connected) Chevalley group $G$ attached to $\Phi$ and $L$ (see 3.3.1) has 
property (CS).}
\item[ii)] {\it The group ${\rm SL}_2$ does not satisfy (CS) (see Corollary 
18 below).}
\end{itemize}

\subsection{ \ }

Let us define {\it projective limits} of groups. Consider a partially 
ordered set $I$ such that, for any $i,j$ in $I$, there is some
$k \in I$ with $k
\geq i$ and $k \geq j$. Assume given a family of groups $G_i$, $i \in I$, 
and morphisms $\varphi_{ji} : G_j \to G_i$, when $j \geq i$, such that
$\varphi_{ii} = {\rm id}$ and $\varphi_{ki} = \varphi_{ji} \circ
\varphi_{kj}$ when $k \geq j \geq i$. By definition, the projective limit
$\underset{i}{\varprojlim} \, G_i$ is the group consisting of families 
$(g_i)_{i \in I}$ such that $g_i \in G_i$ and $\varphi_{ji} (g_j) = g_i$ if 
$j \geq i$.

\medskip

When $\Gamma \subset G({\mathbb Q})$ is an arithmetic group, we can consider 
two projective limits. The first one is
$$
\hat\Gamma = \underset{N}{\varprojlim} \, \Gamma / N \, ,
$$
where $N$ runs over all normal subgroups of finite index in $\Gamma$. We can 
also define
$$
\tilde\Gamma = \underset{a}{\varprojlim} \, \Gamma / \Gamma (a) \, ,
$$
where $\Gamma (a)$, $a \geq 1$, runs over all congruence subgroups of $\Gamma$.
There is a surjective map
$$
\hat\Gamma \to \tilde\Gamma
$$
and we let $C(\Gamma)$ be the kernel of this map. The group $C(\Gamma)$ is 
trivial iff $G$ has property (CS).

Note also that we have an inclusion
$$
\tilde\Gamma \to \widetilde{{\rm GL}_N ({\mathbb Z})} = {\rm GL}_N ({\mathbb 
A}_f) \, ,
$$
where ${\mathbb A}_f = \underset{a}{\varprojlim} \, {\mathbb Z} / a {\mathbb
Z}$ is the ring of finite adeles. In \cite{soule:BMS} (16.1), Bass, Milnor and 
Serre considered the following properties:
\begin{itemize}
\item[a)] $C(\Gamma)$ is finite;
\item[b)] the image of $\hat\Gamma \to {\rm GL}_N ({\mathbb A}_f)$ contains 
a congruence subgroup of ${\rm GL}_N ({\mathbb Z})$.
\end{itemize}

They conjectured that a) and b) are true when $G$ is simple and simply 
connected over ${\mathbb Q}$ (and not necessarily split). Under these 
assumptions, the assertion a) is known today in many cases: see \cite{soule:P} 
\S~9.5, where $G$ can also be defined over some number field. And the 
assertion b) is true when $G({\mathbb R})$ is not compact, by the strong 
approximation theorem (\cite{soule:BMS} loc. cit., \cite{soule:P} Th. 7.12).

\subsection{ \ }

The interest of a) and b) is the following ``rigidity'' result (\cite{soule:BMS} 
Theorem 16.2):

\medskip

\noindent {\bf Proposition 14.} {\it Assume $G$ is a semi-simple group which 
is simply connected ({\rm i.e.} $G$ does not have any nontrivial central 
extension), let $\Gamma \subset G({\mathbb Q})$ be an arithmetic subgroup 
satisfying {\rm a)} and {\rm b)} in {\rm 7.2}, and let
$$
\rho : \Gamma \to {\rm GL}_N ({\mathbb Q})
$$
be any representation. Then there exists an algebraic group morphism
$$
\varphi : G \to {\rm GL}_N
$$
and a subgroup of finite index $\Gamma' \subset \Gamma$ such that the 
restrictions of $\rho$ and $\varphi$ to $\Gamma'$ coincide.}

\medskip

\noindent {\bf Remark.} Stronger results were obtained later by Margulis 
\cite{soule:Ma}; see below Theorem 22.

\subsection{ \ }

We derive from Proposition 14 several consequences.

\medskip

\noindent {\bf Corollary 15.} {\it Let $G$ and $\Gamma$ be as in Proposition 
$14$ and let
$$
\Gamma \to {\rm Aut} (V)
$$
be any representation of $\Gamma$ on a finite dimensional ${\mathbb 
Q}$-vector space. Then $V$ contains a lattice stable by $\Gamma$.}

\medskip

\noindent {\it Proof.} Let $\varphi : G \to {\rm Aut} (V)$ and $\Gamma' 
\subset \Gamma$ be chosen as in the proposition. Then $\varphi (\Gamma')$ is 
contained in an arithmetic subgroup of ${\rm Aut} (V)$ (see Proposition 3), 
hence there is a lattice $\Lambda'$ in $V$ stable by $\Gamma'$ (or a finite 
index subgroup). Let $S \subset \Gamma$ be a set of representatives of 
$\Gamma$ modulo $\Gamma'$. The lattice
$$
\Lambda = \sum_{s \in S} s(\Lambda')
$$
in $V$ is stable by $\Gamma$.

\medskip

\noindent {\bf Corollary 16.} {\it Let $G$ and $\Gamma$ be as in Proposition 
$14$ and let
$$
0 \to V' \to V \to V'' \to 0
$$
be an exact sequence of finite dimensional representations of $\Gamma$ over
${\mathbb Q}$. This sequence splits.}

\medskip

\noindent {\it Proof.} Choose $\Gamma' \subset \Gamma$ such that the 
restriction of the exact sequence to $\Gamma'$ is induced by an exact 
sequence of algebraic representations of $G$. Since $G$ is semi-simple, 
hence reductive, this sequence of representations is split by a section
$\sigma' : V'' \to V$ which commutes with the action of $G$ and $\Gamma'$. 
If $S$ is a set of representatives of $\Gamma$ modulo $\Gamma'$, the formula 
$$
\sigma (x) = \frac{1}{{\rm Card} (S)} \sum_{s \in S} s \, \sigma' \, s^{-1} 
(x) \, ,
$$
$x \in V''$, defines a $\Gamma$-equivariant splitting of the exact sequence.

\medskip

\noindent {\bf Corollary 17.} {\it When $G$ and $\Gamma$ are as in 
Proposition $14$, the abelian group $\Gamma / [\Gamma , \Gamma]$ is finite.}

\medskip

\noindent {\it Proof.} The quotient $\Gamma / [\Gamma , \Gamma]$ of $\Gamma$ 
by its commutator subgroup is abelian and finitely generated. If it was 
infinite there would exist a nontrivial morphism
$$
\chi : \Gamma \to {\mathbb Z} \, .
$$
Let $V = {\mathbb Q}^2$ be equipped with the $\Gamma$-action $\Gamma \to 
{\rm Aut} (V)$ which maps $\gamma$ to $\begin{pmatrix} 1 &\chi (\gamma) \\ 0 
&1 \end{pmatrix}$. We get an exact sequence
$$
0 \to V' \to V \to V'' \to 0
$$
where $\Gamma$ acts trivially on $V' \simeq V'' \simeq {\mathbb Q}$. Since 
$\chi$ is nontrivial, this sequence is not trivial, and this contradicts 
Corollary 16.

\medskip

\noindent {\bf Corollary 18.} {\it The group ${\rm SL}_2$ does not satisfy 
(CS).}

\medskip

\noindent {\it Proof.} Let $\Gamma \subset {\rm SL}_2 ({\mathbb Z})$ be any 
arithmetic subgroup. We shall prove that $C(\Gamma)$ is infinite. If 
$\Gamma' \subset \Gamma$ is a torsion free subgroup 
of finite index, the morphism
$$
C(\Gamma') \to C(\Gamma)
$$
has finite kernel and cokernel, therefore we can assume $\Gamma' = \Gamma$. 
But then, by Theorem 10, $\Gamma$ is free, therefore $\Gamma / [\Gamma , 
\Gamma]$ is a nontrivial free abelian group. The group ${\rm SL}_2$ 
satisfies the strong approximation theorem, therefore b) in 7.2 is true. 
From Proposition 14 and Corollary 17, we conclude that a) is not true, {\it 
i.e.} $C(\Gamma)$ is infinite.

\section{Kazhdan's property (T)}

Let $G$ be a topological group and $\pi$ a unitary representation of $G$ in 
a Hilbert space ${\mathcal H}$. We say that $\pi$ {\it contains almost 
invariant vectors} when, for every $\varepsilon > 0$ and every compact 
subset $K \subset G$, there is a vector $v \in {\mathcal H}$, $v \ne 0$, 
such that
$$
\Vert \pi (g) \, v - v \Vert < \varepsilon
$$
for all $g \in K$.

The group $G$ has {\it property (T)} when any unitary representation $\pi$ 
which contains almost invariant vectors has an invariant vector (a $w \ne 0$ 
such that $\pi (g) \, w = w$ for all $g \in G$).

\medskip

\noindent {\bf Theorem 19} (\cite{soule:Ka}).
\begin{itemize}
\item[i)] {\it Assume that $G$ is locally compact and that $\Gamma \subset 
G$ is a closed subgroup such that the invariant volume of $\Gamma \backslash 
G$ is finite. Then $G$ has property (T) iff $\Gamma$ has property (T).}
\item[ii)] {\it Assume $\Gamma$ is discrete and has property (T). Then 
$\Gamma$ is finitely generated and $\Gamma / [\Gamma , \Gamma]$ is finite.}
\end{itemize}

\medskip

\noindent {\bf Theorem 20} (\cite{soule:OV} Theorem 3.9, p.19). 
{\it Let $G$ be a simple connected
Lie group. Then $G$ has property (T) iff it is not locally isomorphic to 
${\rm SO} (n,1)$ or ${\rm SU} (n,1)$, $n \geq 2$.}

\medskip

(Recall that $G$ is simple if it does not contain any proper nontrivial 
closed normal connected subgroup).

We can combine Theorem 19 i), Corollary 5 ii) and Theorem 20 to show that 
some arithmetic groups have property (T). For instance, ${\rm SL}_N 
({\mathbb Z})$ has property (T) iff $N \geq 3$. For an ``effective'' version 
of that result, see \cite{soule:K}.

\section{Arithmeticity}

\subsection{ \ }

When $G$ is semi-simple over ${\mathbb Q}$, we know from Corollary 5 ii) 
that any arithmetic subgroup $\Gamma \subset G({\mathbb Q})$ has finite 
covolume in $G({\mathbb R})$. A famous conjecture of Selberg asked for a 
converse to this assertion. It was proved by Margulis \cite{soule:Ma1}. We state 
his theorem for simple Lie groups.

\medskip

\noindent {\bf Theorem 21} (\cite{soule:Ma1}). {\it Let $H$ be a connected simple 
non-compact Lie group of rank bigger than one, and $\Gamma \subset H$ a 
discrete subgroup of finite covolume. Then $\Gamma$ is ``arithmetic''.}

\medskip

We need to explain what being ``arithmetic'' means: there exists a linear
algebraic group $G$ over ${\mathbb Q}$, an arithmetic subgroup $\Gamma'$ of 
$G$, a compact Lie group $K$ and an isomorphism of Lie groups
$$
G({\mathbb R}) \simeq H \times K
$$
such that the first projection of $\Gamma'$ into $H$ has finite index in 
$\Gamma$.

\subsection{ \ }

When $H = {\rm PSL}_2 ({\mathbb R})$ (a case of rank one),
Theorem 21 is not true anymore. Indeed,
let $M$ be a compact Riemann surface. Uniformization gives an embedding of 
$\Gamma = \pi_1 (M)$ into ${\rm PSL}_2 ({\mathbb R})$, and the quotient 
$\Gamma \backslash {\rm PSL}_2 ({\mathbb R})$ is compact. But, in general, 
$\Gamma$ is not arithmetic.

\subsection{ \ }

The proof of Theorem 21 uses the following ``superrigidity'' theorem, and 
its non-archimedean analogs (see \cite{soule:Ma2}, \cite{soule:OV} Theorem 6.2.1 or
\cite{soule:Z} for a general statement):

\medskip

\noindent {\bf Theorem 22.} {\it Let $\Gamma \subset H$ be as in Theorem 
$21$. Assume that $G$ is a semi-simple algebraic group over ${\mathbb R}$ 
and $f : \Gamma \to G({\mathbb R})$ is a group morphism such that 
$f(\Gamma)$ is Zariski dense. Then $f$ is the restriction to $\Gamma$ of a 
morphism of Lie groups $H \to G({\mathbb R})$.}

\subsection{ \ }

Let us conclude this survey with another result of 
Margulis \cite{soule:Ma3},
concerning all normal subgroups of a given arithmetic group:

\medskip

\noindent {\bf Theorem 23.} {\it Assume $G$ is a linear algebraic group over 
${\mathbb R}$ such that $G({\mathbb R})$ is connected, simple, not compact 
and of real rank bigger than one. If $\Gamma \subset G({\mathbb R})$ is 
discrete with finite covolume, any normal subgroup $N \subset \Gamma$ has 
finite index in $\Gamma$ or it is contained in the center of $\Gamma$.}

\newpage
\centerline{\Large \bf Appendix}

\bigskip
\bigskip

Following E. Cartan, Cremmer and Julia give in \cite{soule:CJ} the following description of the simple
complex Lie algebra of type $E_7$ and its fundamental representation of dimension 56.

Let $W = {\mathbb C}^8$, with basis $e_i$, $1 \leq i \leq 8$, and $W^*$ its complex dual, 
with dual basis $e_i^*$, $1 \leq i \leq 8$. For any positive integer $k$, we let $\Lambda^k \, W$ 
be the $k$-th exterior power of $W$, i.e. the linear subspace of $W^{\otimes k}$ consisting
 of fully antisymmetric tensors. A basis of $\Lambda^k \, W$ consists of the vectors
$$
e_{i_1} \wedge \ldots \wedge e_{i_k} = \sum_{\sigma \in {\mathcal S}_k} \varepsilon (\sigma)
 \, e_{i_{\sigma (1)}} \otimes \ldots \otimes e_{i_{\sigma (k)}} \, ,
$$
with $1 \leq i_1 < i_2 < \ldots < i_k \leq 8$, where ${\mathcal S}_k$ is the permutation group 
on $k$ letters and $\varepsilon (\sigma)$ is the signature 
of $\sigma$. The exterior product
$$
\Lambda^k \, W \otimes \Lambda^{\ell} \, W \to \Lambda^{k+\ell} \, W
$$
sends $(v_1 \wedge \ldots \wedge v_k) \otimes (w_1 \wedge \ldots \wedge w_{\ell})$ to 
$v_1 \wedge \ldots \wedge v_k \wedge w_1 \wedge \ldots \wedge w_{\ell}$. The basis
 $e_1 \wedge \ldots \wedge e_8$ gives an identification $\Lambda^8 \, W = {\mathbb C}$ and, together
  with the exterior product, an isomorphism
$$
(\Lambda^k \, W)^* = \Lambda^{8-k} \, W
$$
for all $k \leq 8$.

On the other hand, we get a pairing
$$
\Lambda^k \, W \otimes \Lambda^k (W^*) \to {\mathbb C}
$$
by sending $(v_1 \wedge \ldots \wedge v_k) \otimes (\lambda_1 \wedge \ldots \wedge \lambda_k)$
 to the determinant of the $k$ by $k$ matrix $(\lambda_j (v_i))_{1 \leq i , j \leq k}$. 
 This pairing identifies $\Lambda^k (W^*)$ with $(\Lambda^k \, W)^*$.

Let now $V = \Lambda^2 (W^*) \oplus \Lambda^2 (W)$, a complex vector space of dimension 56. 
The complex Lie algebra $\Lambda =  sl_8 ({\mathbb C})$ acts upon $W$, hence on $V$.

Let $\Sigma = \Lambda^4 \, W$, so that $\dim_{\mathbb C} (\Sigma) = 70$. 
From the previous discussion, we get natural pairings
$$
\Lambda^4 \, W \otimes \Lambda^2 (W^*) \overset{\sim}{\longrightarrow}
 (\Lambda^4 \, W)^* \otimes \Lambda^2 (W^*) \to \Lambda^6 (W^*) = \Lambda^2 (W)
$$
and
$$
\Lambda^4 \, W \otimes \Lambda^2 \, W \to \Lambda^6 \, W = (\Lambda^2 \, W)^* = \Lambda^2 (W^*) \, .
$$
Let
$$
\Sigma \otimes V \to V
$$ be the action of $\Sigma$ on $V$ obtained by taking the direct sum of these maps 
and multiplying the result by 2.

The action of
$$
{\mathcal G} = \Lambda \oplus \Sigma
$$
on $V$ defines an embedding
$$
{\mathcal G} \subset {\rm End} (V) \, ,
$$
which is the one given by formulae (B1) in \cite{soule:CJ} (the factor 2 above comes from 
the permutation of $k$ and
$\ell$ in the expression $\sum_{ijk\ell} x^{k\ell}$ of loc.cit.). The vector space ${\mathcal G}$ is stable
 under the Lie bracket, which is given by formulae (B2) in \cite{soule:CJ}, and the action of ${\mathcal G}$ on $V$
  respects the canonical symplectic form on $V$ coming from the pairing
$$
\Lambda^2 (W^*) \otimes \Lambda^2 \, W \to {\mathbb C} \, .
$$
Therefore ${\mathcal G}$ is contained in $sp_{56} ({\mathbb C})$.

Let us now apply Chevalley's construction to this representation of ${\mathcal G}$ on $V$. A Cartan subalgebra
 of ${\mathcal G}$ is the diagonal subalgebra ${\mathcal H} \subset \Lambda$. We let
$$
\varepsilon_i : {\mathcal H} \to {\mathbb C} \ , \quad 1 \leq i \leq 8 \, ,
$$
be the character sending a diagonal matrix to its $i$-th entry. Note that $\varepsilon_1 + 
\varepsilon_2 + \ldots + \varepsilon_8 = 0$ on ${\mathcal H}$. The action of ${\mathcal H}$ on
${\mathcal G} = \Lambda \oplus \Sigma$ is the restriction of the action of $\Lambda =  sl
_8 ({\mathbb C})$.
 Therefore the roots of ${\mathcal H}$ are of two types.

The roots of ``type $\Lambda$'' are those given by the action of ${\mathcal H}$ on $\Lambda$.
These are $\alpha = \varepsilon _i - \varepsilon_j$, for all $i \ne j$, $1 \leq i,j \leq 8$. 
The corresponding eigenspace ${\mathcal G}_{\alpha}$ is spanned by $X_{\alpha} = X_{ij}$,
the matrix having 1 as $(i,j)$ entry, all others being  zero. There are 63 roots of type $\Lambda$.

The roots of ``type $\Sigma$'' are those given by the action of ${\mathcal H}$ on
$\Sigma = \Lambda^4 \, W$. Given four indices $1 \leq i_1 < i_2 < i_3 < i_4 \leq 8$
we get the root $\alpha = \varepsilon_{i_1} + \varepsilon_{i_2} + \varepsilon_{i_3} + e_{i_4}$, with ${\mathcal G}_{\alpha}$ spanned by
$$
X_{\alpha} = \frac{1}{2} \, e_{i_1} \wedge e_{i_2} \wedge e_{i_3} \wedge e_{i_4} \, .
$$
There are 70 roots of type $\Sigma$. Let $\Phi$ be the set of all roots.

We claim that the vectors $X_{\alpha}$ and $H_{\alpha} = [X_{\alpha} , X_{-\alpha}]$,
 $\alpha \in \Phi$, form a Chevalley basis of ${\mathcal G}$. According to \cite{soule:H}, 
 proof of Proposition 25.2, this will follow if we prove that the Cartan involution $\sigma$ satisfies
\begin{equation}
\label{soule:eq1new}
\sigma (X_{\alpha}) = -X_{-\alpha}
\end{equation}
and that the Killing form $K$ is such that
\begin{equation}
\label{soule:eq2new}
K (X_{\alpha} , X_{-\alpha}) = 2/(\alpha , \alpha) \, ,
\end{equation}
for every root $\alpha \in \Phi$.

The Cartan involution $\sigma$ on ${\mathcal G}$ is the restriction of the Cartan 
involution on ${\rm End} (V)$, so it is the standard one on $\Lambda =  sl_8 ({\mathbb C})$ and we get
$$
\sigma (X_{ij}) = - X_{ji} \, .
$$
On the other hand, the pairings
 of $Z = \Lambda^4 \, W$ with $\Lambda^2 (W^*)$ and $\Lambda^2 \,
 W$ are dual to each other. Therefore, if $x \in Z$ we have $\sigma (x) = -x^*$, where $x^*$ is the image of $x$ by the isomorphism
$$
\Lambda^4 \, W \overset{\sim}{\longrightarrow} (\Lambda^4 \, W)^* = \Lambda^4 (W^*)
$$
followed by the identification of $W$ and $W^*$ coming from the chosen bases. This sends 
$e_1 \wedge e_2 \wedge e_3 \wedge e_4$ to $e_5^* \wedge e_6^* \wedge e_7^* \wedge e_8^*$, 
and then to $e_5 \wedge e_6 \wedge e_7 \wedge e_8$. We conclude that
$$
\sigma (e_1 \wedge e_2 \wedge e_3 \wedge e_4) = - e_5 \wedge e_6 \wedge e_7 \wedge e_8 \, .
$$
The root corresponding to $e_5 \wedge e_6 \wedge e_7 \wedge e_8$ is
$$
\varepsilon_5 + \varepsilon_6 + \varepsilon_7 + \varepsilon_8 = - 
(\varepsilon_1 + \varepsilon_2 + \varepsilon_3 + \varepsilon_4) \, .
$$
This proves (\ref{soule:eq1new}) for roots of type $\Sigma$.

Let us now check (\ref{soule:eq2new}). According to the definitions in \cite{soule:H} 8.2 and 8.3, we have
$$
(\alpha , \alpha) = K (T_{\alpha} , T_{\alpha})
$$
where $T_{\alpha} \in {\mathcal H}$ is defined by the equality
$$
\alpha (h) = K(T_{\alpha} , H)
$$
for all $H \in {\mathcal H}$. When $X,Y \in \Lambda$ are two $8 \times 8$ matrices of trace zero,
we have, as indicated in \cite{soule:CJ} (B5),
$$
K(X,Y) = 12 \, {\rm tr} (XY) \, .
$$
Let $\alpha = \varepsilon_i - \varepsilon_j$ be a root of type $\Lambda$ and $H_{ij}$ the
 diagonal matrix such that $\varepsilon_i (H_{ij}) = 1$, $\varepsilon_j (H_{ij}) = -1$ and
 $\varepsilon_k (H_{ij}) = 0$ if $k \notin \{ i,j \}$. For any $H \in {\mathcal H}$ we have
$$
\alpha (H) = {\rm tr} \, (H_{ij} \, H) \, ,
$$
therefore
$$
T_{\alpha} = H_{ij} / 12
$$
and
$$
(\alpha , \alpha) = \frac{1}{144} \, K (H_{ij} , H_{ij}) = \frac{24}{144} = \frac{1}{6} \, .
$$
Since $K (X_{\alpha} , X_{-\alpha}) = 12$, the equality (\ref{soule:eq2new}) holds true.

Assume now that $\alpha = \varepsilon_1 + \varepsilon_2 + \varepsilon_3 + \varepsilon_4$,
 $X_{\alpha} = \frac{1}{2} \, e_1 \wedge e_2 \wedge e_3 \wedge e_4$ and
 $X_{-\alpha} = \frac{1}{2} \, e_5 \wedge e_6 \wedge e_7 \wedge e_8$.
  According to (B5) in \cite{soule:CJ} we have
$$
K (X_{\alpha} , X_{-\alpha}) = \frac{2}{24} \, \frac{1}{4} \, (4!) (4!) = 12 \, .
$$
Let $H'$ be the diagonal matrix such that $\varepsilon_i (H') = 1/2$ when $1 \leq i \leq 4$
and $\varepsilon_i (H') = -1/2$ when $5 \leq i \leq 8$. Given any $H$ in ${\mathcal H}$ we have
$$
{\rm tr} \, (H'H) = \frac{1}{2} (\varepsilon_1 + \varepsilon_2 + \varepsilon_3 + \varepsilon_4 -
 \varepsilon_5 - \varepsilon_6 - \varepsilon_7 - \varepsilon_8) (H) = \alpha (H) \, .
$$
Therefore $T_{\alpha} = H' / 12$ and
$$
(\alpha , \alpha) = \frac{1}{144} \, K(H',H') = \frac{12}{144} \times \frac{8}{4} = \frac{1}{6} \, .
$$
Therefore (\ref{soule:eq2new}) is true for $\alpha$.

The Lie algebra ${\mathcal G}$ is simple of type $E_7$. Indeed, a basis of its root system $\Phi$
 consist of $\alpha_1 = \varepsilon_1 - \varepsilon_2$, $\alpha_2 = \varepsilon_4 + \varepsilon_5 + 
 \varepsilon_6 + \varepsilon_7$, $\alpha_3 = \varepsilon_2 - \varepsilon_3$, $\alpha_4 = \varepsilon_3 - 
 \varepsilon_4$, $\alpha_5 = \varepsilon_4 - \varepsilon_5$, $\alpha_6 = \varepsilon_5 - \varepsilon_6$ and 
 $\alpha_7 = \varepsilon_6 - \varepsilon_7$. Its Dynkin diagram is the one of $E_7$ (\cite{soule:H}, 11.4).

Let us now consider the representation $\rho$ of ${\mathcal G}$ on $V$.
 Its weight vectors are
$e_i^* \wedge e_j^* \in \Lambda^2 (W^*)$ and $e_i \wedge e_j \in \Lambda^2 \, W$, $1 \leq i < j \leq 8$,
 with corresponding weights $-\varepsilon_i - \varepsilon_j$ and $\varepsilon_i + \varepsilon_j$. 
 The root lattice $L_0$ of $E_7$ has index 2 in its weight lattice $L_1$ (\cite{soule:H}, 13.1). Since
  the weights of $\rho$ are not in $L_0$ they must span the lattice $L_1$. Therefore, the
   Chevalley group $G$ generated by the endomorphisms $\exp (t\rho (X_{\alpha}))$,
   $t \in {\mathbb C}$, $\alpha \in \Phi$, is the simply connected Chevalley group of type $E_7$.
   Its set of real points $G({\mathbb R})$ is the real Lie group spanned by the endomorphisms
    $\exp (t\rho (X_{\alpha}))$, $t \in {\mathbb R}$, $\alpha \in \Phi$ (\cite{soule:St}, \S~5, Th.~7, Cor.~3),
    i.e. the split Lie group $E_{7(+7)}$.

Let $M \subset V$ be the standard lattice, with basis $e_i^* \wedge e_j^*$ and $e_i \wedge e_j$,
 $1 \leq i < j \leq 8$. The group $E_7 ({\mathbb Z}) = E_{7(+7)} \cap {\rm Sp}_{56} ({\mathbb Z})$
  is the stabilizer of $M$ in $G$. So, according to \cite{soule:St}, \S~8, Th.~18, Cor.~3, to check that
   $E_7 ({\mathbb Z}) = G ({\mathbb Z})$, all we need to prove is that the lattice $M$ is
    admissible, i.e. stable by the endomorphisms $\rho (X_{\alpha})^n / n!$ for all $n \geq 1$ and $\alpha \in \Phi$.

When $\alpha = \varepsilon_i - \varepsilon_j$ is of type $\Lambda$, $\rho (X_{\alpha}) = X_{ij}$ has
square zero and stabilizes the standard lattice $M$. Assume finally that
 $\alpha = \varepsilon_1 + \varepsilon_2 + \varepsilon_3 + \varepsilon_4$, hence 
 $X_{\alpha} = \frac{1}{2} \, e_1 \wedge e_2 \wedge e_3 \wedge e_4$. By definition
  of the action of $Z$ on $V = \Lambda^2 (W^*) \oplus \Lambda^2 \, W$, 
  $\rho (X_{\alpha})$ sends $e_i \wedge e_j$ to $\pm \, e_k^* \wedge e_{\ell}^*$ when
  $5 \leq i < j$, $k < \ell$ and $\{ i,j,k,\ell \} = \{5,6,7,8\}$. When $i < 5$,
   $\rho (X_{\alpha}) (e_i \wedge e_j) = 0$. Similarly, when $i < j \leq 4$,
   $\rho (X_{\alpha})$ sends $e_i^* \wedge e_j^*$ to $\pm \, e_k \wedge e_{\ell}$
   with $\{ i,j,k,\ell \} = \{ 1,2,3,4 \}$, and $\rho (X_{\alpha}) (e_i^* \wedge e_j^*) = 0$
    if $j > 4$. From this it follows that the endomorphism $\rho (X_{\alpha})$ has
    square zero and stabilizes $M$. Therefore $E_7 ({\mathbb Z}) = G({\mathbb Z})$.

\end{document}